\documentclass[11pt,a4paper,reqno]{amsproc}

\usepackage{graphicx}
\usepackage{amsfonts}
\usepackage{amsmath,amscd}
\usepackage{amsthm}

{ \theoremstyle{plain}
\newtheorem{theorem}{Theorem}
\newtheorem{lemma}{Lemma}

}

{ \theoremstyle{definition}
}

{ \theoremstyle{definition}

}

{ \theoremstyle{remark}

}

\def\scal#1#2{\langle #1, #2\rangle }

\def\R{{\mathbb R}}
\def\C{{\mathbb C}}
\def\Z{{\mathbb Z}}

\def\prx{\frac{\partial }{\partial x}}
\def\pry{\frac{\partial }{\partial y}}

\DeclareMathOperator{\im}{\mathrm{Im}}
\DeclareMathOperator{\re}{\mathrm{Re}}
\DeclareMathOperator{\sn}{\mathrm{sn}}
\DeclareMathOperator{\sh}{\mathrm{sinh}}
\DeclareMathOperator{\ch}{\mathrm{cosh}}

\begin{document}

\title[Entire one-periodic maximal surfaces]{Entire one-periodic maximal surfaces
\footnote{\textbf{Printed in}: \textit{The Mansfield-Volgograd anthology}, Ed.: James York Glimm, Mansfield University of Pennsylvania, 2000. 148-156.}}

%
%

\author{Sergienko Vladimir V.}

\email{}
\thanks{The first author was supported by grant INTAS No.10170}

%
%
\author{Tkachev Vladimir G.}


\thanks{The second author was supported by grant INTAS No.10170 and
Ministerstvo Vyshego Obrazovaniya Rossii N 97-0-1.3-114}

%
%

%
%

\begin{abstract}
In the present paper we study two-dimensional
maximal surfaces with harmonic level-sets. As a corollary we
obtain a new class of one-periodic maximal surfaces.
\end{abstract}

\maketitle

\section{Introduction}

Let $\R _1^{n+1}$ be $(n+1)$-dimensional Minkowski space with the
standard metric
$$
\scal{\chi '}{\chi ''}=x'_1\cdot x''_1+\ldots +x'_n\cdot x''_n-t'\cdot t'',
$$
where $\chi =(x_1,\ldots , x_n, t)$ is a point in $\R _1^{n+1}$.

Let $M$ be a surface in $\R _1^{n+1}$ given by
$t=f(x_1,\ldots ,x_n)$, and $f$ be a  function in a domain
$\Omega \subset \R ^{n}\equiv \{\chi : t=0\}$. We shall assume that
$f$ is $C^2$-smooth everywhere in $\Omega$ except for a set $A\subset
\Omega$ consisting of isolated points only.

A surface $M$ is {\it space-like} if the induced from $\R_1^{n+1}$
metric is the positive definite metric. This is equivalent to the
following inequality
\begin{equation}
|\nabla f|^2=\sum_{i=1}^{n}\,f'^2_{x_i}(x)<1, \quad \forall x\in
\Omega \setminus A,
\label{second}
\end{equation}
where the lower index denotes a partial derivative
with respect to the corresponding variable: $f_{x_i}=\frac{\partial f}{\partial x_i}$,
and $|\nabla f|^2=\sum_{i=1}^{n}{f}'^2_{x_i} $.

A space-like surface $M$ is called {\it maximal} if
the following equality holds
\begin{equation}
\sum_{i,j=1}^{n}\,f''_{ij}\, \biggl(\delta _{ij} (1-|\nabla
f|^2)+f'_{x_i}f'_{x_j}\biggr)=0
\label{first}
\end{equation}
everywhere in $\Omega \setminus A$.
Then it is well-known  that (\ref{first}) is equivalent to the
vanishing of the mean curvature of $M$ (with respect to
its embedding in the Minkowski space  $\R_1^{n+1}$).

Cheng and Yau in \cite{CY2} proved that for every entire maximal
surface $M$ (i.e. the surface to be defined over the whole $\R^{n}$)
satisfies the Bernstein property. In other words, an entire solution $f(x)$
to~(\ref{second})-(\ref{first}) is always an affine function.
On the other hand, the study of {\it almost-entire} solutions to
(\ref{second})-(\ref{first}), i.e.  solutions that they are
of $C^2$ outside of a non-empty descrete set $A$, is of great interest. Really, there are a lot of connections between
this theory and the modern physics (see \cite{BL}).

The first break-through in this direction was due to Ecker
\cite{Ec}, who established that the rotationally symmetric maximal
surfaces (`maximal catenoids')
\begin{equation}
\|x\|=\sqrt{x_1^2+\ldots +x_n^2}, \quad
t=c\int\limits_{0}^{\|x\|}\left( c^2+\lambda
^{2(n-1)}\right)^{-\frac{1}{2}}\; d\lambda,
\label{third}
\end{equation}
are only almost-entire solutions with $A$ consisting of the
origin. Moreover, in the same paper it was proved that every isolated
singular point of a maximal (not necessarily entire) surface
behaves as a light cone. Namely, let $a\in \R^{n}$ be an isolated
singular point of a solution $f$ which is defined at a
neighborhood of $a$ then $f$ can be extended continuously at $a$ and
$$
f(x)-f(a)=\pm \|x-a\|+o(\|x-a\|).
$$
An important characteristic of a singular point is the flux of the
solution which is defined as
\begin{equation}
\mu _f(a)=\int\limits _{C_{a}}\frac{\scal{\nabla f}{\nu}}{\sqrt{1-|\nabla
f|^2}},
\label{fourth}
\end{equation}
where the integral is taken over a closed surface $C_a\subset \R^{n}$ which encloses $a$ and contains no singular points,
and $\nu$ denotes the outward unit normal to $C_a$. It follows
from (\ref{first}) that integral (\ref{fourth}) is independent on a choice of $C_a$.
It was shown in \cite{Mik} that $\mu_f(a)$ is a Lorentzian invariant of
$M$. It is not hard to prove that a point $a\in \Omega$ is an
essential singularity of the solution $f$ if and only if
$\mu_f(a)\ne 0$ \cite{M3}.

In recent papers \cite{KaM1}, \cite{Ka} the asymptotic behaviour
and the existence questions for solutions to
(\ref{second})-(\ref{first}) were studied. In particularly,
it was shown in \cite{KaM1} that under some natural geometrical assumptions on
the finite set $A$ there exists a unique almost-entire
solution $f$ with the prescribed fluxes $\mu_f(a_i)$, $a_i\in A$, where
$A$ is the singular set of the solution.

On the other hand, there are no explicit examples except for
the mentioned above "maximal catenoids" (\ref{third}) even in the
two-dimensional case. In this paper we construct a one-parametric
family of periodic almost-entire maximal surfaces whose
singular set is discrete (consists of isolated points)
and located on the fixed line. These
examples are based on the following general assertion which
completely characterizes all maximal surfaces with
harmonic level-sets.

\begin{theorem}
Let $f(x,y)$ satisfy {\rm (\ref{second})-(\ref{first})} and let
$f(x,y)=F(\varphi(x,y))$, where $\varphi (x,y)$ is a harmonic
function. Then $\varphi (x,y)$ is a real part of the holomorphic
function $h(w)=\int \frac{dw}{g(w)}$, where $g(w)$ is one of the
following

(i) $g(w)=aw+c$;

(ii) $g(w)=ae^{bw}$;

(iii) $g(w)=a\sin(bw+c)$.\\
Here $w=x+iy$ and $a^2$, $b^2 \in \R$, $c\in {\bf C}$.
\end{theorem}

The cases (i) and (ii) lead us to the well-known examples:
the plane, the rotational surface, the helicoid and
Sherk's maximal surfaces. In the case (iii) the surface is
space-like only if $a,b \in \R$, moreover in that case it has infinitely many isolated
singularity points. The explicit expression and analysis we shall give
in the remained part of the paper.

In fact, the examples of one-periodic maximal surfaces
constructed in this paper are a small part of a bigger family of
double- and one-periodic maximal two-dimensional surfaces which we treat in the forthcoming
paper ~\cite{ST2}.

\section{Preliminaries}

Let $u(x,y)$ be a $C^2$-function such that  $|\nabla
u(x,y)|<1$ and
\begin{equation}
(1-{u'}^2_y){u''}_{xx}+2{u'}_x{u'}_y{u''}_{xy}+(1-{u'}^2_x){u''}_{yy}=0.
\label{ur1}
\end{equation}
Then the graph of $u(x,y)$ is a maximal surface in $\R_1^3$.

Now we consider the function $u(x,y)$ such that $u(x,y)=F(\varphi (x,y))$,
where $F=F(\eta)$ and $\varphi(x,y)$~ are some twice-differentiable functions. Then
(\ref{ur1}) can be brought to the following form:
\begin{equation}
A(x,y)F''_{\eta  \eta }+B(x,y)F'_{\eta }+C(x,y){F'_{\eta
}}^3=0,
\label{ur2}
\end{equation}
where $A(x,y)={\varphi' }^2_x+{\varphi' }^2_y$, $B(x,y)={\varphi''
}_{xx}+{\varphi'' }_{yy}$, $C(x,y)=-{\varphi' }^2_x{\varphi''
}_{yy}+2{\varphi' }_x{\varphi' }_y{\varphi'' }_{xy}-{\varphi'
}^2_y{\varphi'' }_{xx}$.

\begin{lemma}
Let $f(w)$, $g(w)$ be holomorphic functions, $w=x+iy$. Then the
following identities take place
$$
\prx \re (f \bar g)=\re (f'\bar g +f \bar g');
$$
\begin{equation}
\pry \re (f \bar g)=-\im (f'\bar g -f \bar g').
\label{l1}
\end{equation}
Here $\bar g$ denotes the conjugate to $g$ function.
\label{lemm1}
\end{lemma}

\begin{proof} The Cauchi-Riemann conditions imply
$$
\prx \re f=\re f', \qquad \pry \re f=-\im f',
$$
\begin{equation}
\prx \im f=\im f', \qquad
\pry \im f=\re f'.
\label{koshi}
\end{equation}

We prove the validity of the first equality only:
\begin{equation*}
\begin{split}
\prx \re (f \bar g)&=\prx (\re f \re g+\im f \im g)
=\re f' \re g+\re f \re g'\\
&+\im f' \im g+\im f \im g'= \re (f'\bar g +f \bar g')
\end{split}
\end{equation*}
The second equality can be proved by the same way, hence the lemma
is proved completely.
\end{proof}

Let $\varphi (x,y)=\re h(w)$, $w=x+iy$, where $h(w)\in {\bf H}(D)$
is a holomorphic function in the domain $D$.
In order to find the coefficients $A$, $B$ and $C$ of equation (\ref{ur2}),
from (\ref{koshi})  we notice
$$
\varphi' _x=\prx \re h=\re h', \qquad
\varphi'' _{xx}=\prx \re h'=\re h'',
$$
$$
\varphi' _y=\pry \re h=-\im h',\qquad
\varphi'' _{xy}=\pry \re h'=-\im h'',\quad
\varphi'' _{yy}=-\re h''.
$$
Then
\begin{equation*}
\begin{split}
A(x,y)&={\varphi }^2_x+{\varphi }^2_y=|h'(w)|^2,
\\
B(x,y)&={\varphi }_{xx}+{\varphi }_{yy}=0,
\\
C(x,y)&=\re (h''{\bar {h'}}^2)
\end{split}
\end{equation*}
and the equation~(\ref{ur2}) becomes
$$
|h'|^2F''+\re (h''
\bar{h'}^2){F'}^3=0.
$$
Setting
\begin{equation}\label{recall}
g(w)\equiv 1/h'(w)
\end{equation}
we find that
\begin{equation}
\frac{F''(\varphi)}{F'^3(\varphi)}=\frac{\re g'}{|g|^2}.
\label{chast1}
\end{equation}

\begin{lemma}
The term $\frac{1}{|g|^2}\,\re g'$ in the equation {\rm
(\ref{chast1})} depends only on $\varphi (x,y)=\re h(w)$  if and only if
\begin{equation}
gg''-{g'}^2=c,
\label{mmm}
\end{equation}
where $c$ is a real constant.
\label{lemm2}
\end{lemma}

\begin{proof} Let $\psi (x,y)=\re (h' \bar h'
g')\equiv \frac{1}{|g|^2}\,\re g'$. We show that the condition of
the functional dependence $\frac{\partial (\varphi (x,y),\psi
(x,y))}{\partial (x,y)}=0$ is equivalent to $gg''-{g'}^2=c$, $c\in
\R{}$. Indeed, by virtue of (\ref{koshi}) $\varphi' _x=\re h'$, $\varphi'
_y=-\im h'$ and (\ref{l1}), we have
\begin{equation}
\begin{split}
\psi' _x &=2\re
(h''\bar h')\re g'+\re (h'\bar h'g''),\\
\psi' _y &=-2\im
(h''\bar  h')\re g'-\im (h'\bar h'g'').
\label{kos1}
\end{split}
\end{equation}

Then
$$
0=\frac{\partial (\varphi ,\psi)}{\partial (x,y)}=\varphi'
_x\psi' _y-\varphi' _y\psi' _x=-\left(2\im (h'' \bar  h') \re g'+\right.
$$
$$
\left.+\im (h'\bar  h'g'')\right)\cdot \re h'+\left(2\re (h'' \bar  h') \re g'
+\re (h'\bar  h'g'')\right)\cdot \im h'=
$$
$$
=-2\im (h'' \bar{h'}^2)\re g'-\im (g''h'\bar{h'}^2),
$$
and, finally,
$$
-2|h'|^4\im \left( \frac{h''}{{h'}^2}\right)\re g'-|h'|^4\im \left(
g''\frac{1}{h'}\right)=0.
$$

Simplifying the last expression yields
$$
0=-2\im g'\re g' +\im (g'' g)=\im (g''g-{g'}^2),
$$
The latter identity holds in a non-empty domain $D$, hence
by the uniqueness theorem for analitic functions, there exists a
real constant $c$ such that $gg''-{g'}^2=c$. The lemma is proved.
\end{proof}

\section{The construction of  examples}

Now we consider the differential equation (\ref{mmm}) with a holomorphic in some domain $D$ function $g(w)$.
One can easily show that the set of solutions of this equation
makes up the following functional family:
(a) $g(w)=aw+c$; (b) $g(w)=ae^{bw}$; (c) $g(w)=a\sin (bw+c)$, $a^2$,
$b^2 \in \R$, $c \in \C$.

The cases (a)-(b) lead us to the classic examples of the maximal surfaces
such as the plane, the rotational surface, the helicoid and
Sherk's maximal surfaces.

Now we consider the last case, when $g(w)=\sin w$\footnote{The general case (c) is reduced to this equation by a suitable isometry and homothety}. Here we have
$h'(w)=\frac{1}{g(w)}$ and, hence,
$$
h(w)=\frac{1}{2}\ln \frac{\cos w -1}{\cos w+1}+const.
$$
Without loss of generality, we may assume that the constant in the
right hand side of the last equality is identically zero. Then
$$
\varphi (x,y)= \re h(w)=\frac{1}{2}\ln \left
|\frac {\cos w -1}{\cos w+1} \right |=\frac{1}{2}\ln\frac{\ch y-\cos
x}{\ch y+\cos x},
$$
and
$$
\frac{1}{|g|^2}\re g'=\frac{\re \cos w }{|\sin w|^2}.
$$
On the other hand,
$$
2\sh 2\varphi (x,y)=\left
|\frac{\cos w -1}{\cos w+1} \right |-\left |\frac{\cos w +1}{\cos
w-1}\right | =
-\frac {4\re \cos w}{|\sin w|^2}=-4\frac{1}{|g|^2}\re g'.
$$
Then the equation (\ref{chast1}) takes the form
$$
F''_{\eta \eta}+\frac{1}{2}{F'_{\eta }}^3\sh 2 \eta =0.
$$
By solving the ordinary differential equation we arrive at
$$
\frac{1}{F'^2_\eta}=\frac{1}{2}\ch 2\eta+\frac{k}{2}, \quad k
\equiv const,
$$
and
$$
F'_\eta(\eta)=\frac{1}{\sqrt{\frac{1}{2}\ch 2\eta
+\frac{k}{2}}}=\frac{1}{\sqrt{\frac{1}{4}(e^{2\eta }+e^{-2\eta})
+\frac{k}{2}}}=\frac{2e^{\eta }}{\sqrt{e^{4\eta }+2ke^{2\eta }+1}}.
$$

To find the admissible values of the parameter $k$ which correspond
to the space-like examples, we check when the inequality
$|\nabla F(\varphi (x,y))|<1$ holds. For this purpose we write
\begin{equation}
|\nabla F(\varphi (w))|=|F'(\varphi)|\,|\nabla \varphi|=
|F'(\varphi)|\,|h'(w)|
\label{ggg}
\end{equation}
and by using a new variable $\gamma =\frac{\cos x}{\ch y}$,
we obtain
$$
|h'(w)|=\frac{1}{|\sin w|}=\frac{1}{\sqrt{\ch^2y-\cos^2x}}=
\frac{1}{\ch y\sqrt{1-\gamma^2}}.
$$

On the other hand, using the exact form of $\varphi$ given above, we find
$$
|F'(\varphi)|=\sqrt{\frac{2(1-\gamma^2)}{(1+k)-(k-1)\gamma^2}}.
$$
Substituting the above expression in (\ref{ggg}) yields
$$
|\nabla F(\varphi (w))|=\frac{1}{\ch
y}\sqrt{\frac{2}{(1+k)-(k-1)\gamma^2}}.
$$

Taking into account, that $\ch y$ and $\gamma=\cos x/\cosh x$ may
change by independent manner, we obtain that the space-likeness
condition takes place only when $k>1$.

We denote $\xi =e^{\eta }$ and assume without loss of generality that $F(0)=0$. Then
$$
F(\eta)=2\int\limits _{1}^{e^{\varphi}} \frac{d\,
\xi}{\sqrt{{\xi}^4 +2k{\xi}^2+1}}=
\sqrt{2}\int\limits _{0}^{\tanh  \eta} \frac{d\,
t}{\sqrt{1-{t}^2}\sqrt{(1+k)-(k-1){t}^2}}=
$$
$$
=\frac{\sqrt{2}}{\sqrt{1+k}}\int\limits _{0}^{\tanh  \eta} \frac{d\,
t}{\sqrt{1-{t}^2}\sqrt{1-{\alpha}^2{t}^2}},
$$
where $\eta=\frac{1-{\xi}^2}{1+{\xi}^2}$ and ${\alpha}^2=\frac{k-1}{k+1}$.
Let us introduce  ${\alpha '}^2=1-{\alpha }^2=\frac{2}{k+1}$.
Then
$$
F(\eta)=\alpha '\int\limits _{0}^{\tanh  \eta} \frac{d\,
t}{\sqrt{1-{t}^2}\sqrt{1-{\alpha}^2{t}^2}},
$$
i.e. by means of the  Jacobi elliptic sinus,  we find
$$
\sn \left( \frac{F(\eta)}{\alpha '};\alpha \right)=\tanh  \eta
$$

\begin{center}
\begin{figure}[h]
\includegraphics[height=0.25\textheight,keepaspectratio=true]{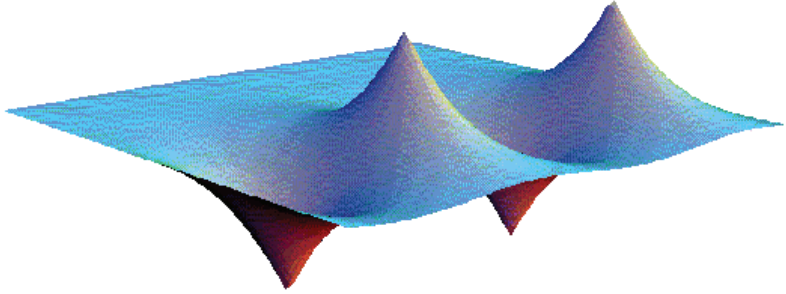}
\caption{One-periodic surface, $\alpha=0.6$} \label{fig:oneper}
\end{figure}
\end{center}

Thus, we have the solution $z=F(\varphi(x,y))$ given by
$$
\sn \left( \frac{z}{\alpha '};\alpha \right)=\tanh  \varphi(x,y)=
-\frac{\cos x}{\ch y}.
$$

The above can be we can summarized as follows.

\begin{theorem}
Let $\alpha \in (0;1)$ and $\alpha '=\sqrt{1-{\alpha}^2}$. Then
the surface $M(\alpha)$ given implicitly by
$$
\sn \left( \frac{z}{\alpha '};\alpha \right)=\frac{\cos x}{\ch y},
$$
is a  maximal surface in ${\bf R}^3$. Moreover, this surface is a graph of a real
analytic function everywhere except for the  set consisting of singular points
$$
A_k=(\pi k;0), \qquad k\in \Z.
$$
\end{theorem}

We observe that for different values of $\alpha\in(0,1)$, the surfaces
$M(\alpha)$ are Lorentz non-isometric. One can also see that $M(\alpha)$ is located in the parallel slab
$|z|\leq \mathrm{K}(\alpha)\alpha'$, where $\mathrm{K}(\alpha)$ is the complete elliptic integral of the first kind (the least positive solution of equation $\sn(\mathrm{K}(\alpha),\alpha)=1$). Moreover, one can show that the flux $\mu(A_k)$ at the singular point $A_k=(\pi k,0)$
is equal to
$$
\mu(A_k)=4\int_{0}^{\pi/2}\frac{\alpha' dt}{\sqrt{1-\alpha'^2\cos^2t}}=
4\alpha'\mathrm{K}(\alpha').
$$

\end{document}